\documentclass{article}
\usepackage[utf8]{inputenc}
\usepackage{amssymb,amsmath,latexsym}
\usepackage[english]{babel}
\usepackage{amscd,amsxtra}
\usepackage{amsthm}
\usepackage{verbatim}
\usepackage{graphicx}
\usepackage{enumitem}
\usepackage{color}
\usepackage{tikz}
\usetikzlibrary{calc, intersections, through, backgrounds}
\usetikzlibrary{shapes.geometric}
\tikzstyle{point}=[ball color=white, circle, draw=black, inner sep=0.1cm]
\tikzstyle{red}=[ball color=red, circle, draw=black, inner sep=0.1cm]
\tikzstyle{blue}=[ball color=blue, circle, draw=black, inner sep=0.1cm]
\usepackage{caption}
\usepackage[rightcaption]{sidecap}
\usepackage[colorinlistoftodos]{todonotes}
\usetikzlibrary{arrows.meta}
\tikzstyle{point}=[ball color=white, circle, draw=black, inner sep=0.1cm]
\usepackage{authblk}
\usepackage{soul}
\usepackage{float}

\newcommand{\verteq}{\rotatebox{90}{$\,=$}}
\newcommand{\equalto}[2]{\underset{\overset{\verteq}{#2}}{#1}}
\newcommand{\equaltog}[2]{\overset{\underset{\verteq}{#2}}{#1}}

\newtheorem{theorem}{Theorem}[section]

\newtheorem{corollary}[theorem]{Corollary}

\newtheorem{definition}[theorem]{Definition}
\newtheorem{example}[theorem]{Example}

\newtheorem{proposition}[theorem]{Proposition}
\newtheorem{remark}[theorem]{Remark}

\DeclareMathOperator{\rad}{rad}

\title{The relation between different edge spans of a graph}
\author[1]{Aljoša Šubašić \thanks{Corresponding author: aljsub@pmfst.hr, ORCID:  0000-0002-6943-0856, Rudjera Bo\v skovi\' ca 33, 21000 Split, Croatia}}
\author[1]{Tanja Vojkovi\'c}
\affil[1]{Faculty of Science, University of Split}

\begin{document}
\maketitle
\begin{abstract}
   In several recent papers, the maximal safety distance that two players can maintain while moving through a graph has been defined and studied using three different spans of the graph, each with different movement conditions. Mainly, vertex spans have been studied, in which players visit all the vertices of a graph. In this paper, we analyze the values of three edge spans, which represent the maximal safety distance that two players can maintain while visiting all the edges of a graph. We present edge span values for some graph classes and examine the relationship between different variants of edge spans. 
\end{abstract}

Keywords: safety distance, graph spans, strong span, direct span, Cartesian span \\
AMS Subject Classification: 05C12, 	05C90

\section{Introduction}
Since the early days of graph theory, various types of walks, trails, and paths have been defined, analyzed, and applied to graphs and networks. Prominent examples include Eulerian circuits, Hamiltonian cycles, Traveling Salesman, and Chinese Postman problems. More contemporary research of walks and paths in a graph include for example \cite{about,aziz} and are often focused on optimizing transport and trade routes \cite{ambrosino2,rumi}. With the emergence of the safety distance concept, in response to the pandemic, many research papers focused on exploring distances in graphs, from the game theory point of view \cite{dash,lagos}, but also interdisciplinary \cite{patasius,rupa}. Recently, Banič and Taranenko examined the situation of two players walking in a graph and analyzed the maximal safety distance they can maintain while visiting either all the vertices or all the edges of the graph \cite{banic}. They approached this idea by introducing three distinct types of vertex and edge spans, where the term 'span' originates from topology \cite{lelek}. Their definitions are based on well-known graph products, strong, direct and Cartesian products, so the span names are inherited. In our previous work \cite{nas}, we approached the definition of graph spans differently than the original authors, directly relating it to the maximal safety distance that two players can maintain while traversing a graph. After demonstrating that our alternative definition yields the same value, we presented results regarding vertex spans of different graphs and explored their interrelations. Since then, several papers have been written about spans, exploring vertex spans of multilayered graphs, \cite{multil}, graphs with spans $1$ and shortest optimal walks, \cite{taranenko1} and strong vertex span of trees, \cite{taranenko2}.\\
This paper serves as a continuation of our previous work, focusing on edge spans. We define edge spans in an analogous manner to the vertex spans described in \cite{nas} and present results for edge spans in various classes of graphs. Our main results consider the relation between different variants of edge spans. We prove that the values of direct and Cartesian edge spans differ by at most one (a result analogous to the one for vertex spans in \cite{nas}), and proceed to prove that the strong edge span is bounded by the values of direct and Cartesian edge spans in a way that it is either equal to the larger of those values, or greater by one. This strong result has also been proven to be true for vertex spans. Moreover, we find a family of graphs for which the upper bound is tight. 

We begin in Section \ref{def} by providing the necessary definitions for approaching edge spans, following an analogous approach to the vertex spans discussed in \cite{nas}. Subsequently, in Section \ref{edge}, we present the three distinct edge spans. Section \ref{prelim} contains preliminary results that will be utilized in the proofs throughout the remainder of the paper. The primary findings on edge spans are presented in Section \ref{res_edge} and the paper ends with a conclusion and further work in Section \ref{con}.

\section{Definitions and previous results}\label{def}

To provide clarity in our proofs, we will continue to refer to two players, Alice and Bob, who traverse a graph by moving through its vertices via its edges. The distance between Alice and Bob is determined by the distance of the vertices they occupy at each step. We will explore three distinct modes of movement for Alice and Bob:
\begin{itemize}
	\item Traditional movement rules: In this mode, during each step, Alice and Bob independently choose to either move to an adjacent vertex or remain in their current position.
	\item Active movement rules: In this mode, during each step, both Alice and Bob are required to move to adjacent vertices simultaneously.
	\item Lazy movement rules: In this mode, during each step, one player moves to an adjacent vertex while the other player remains stationary.
\end{itemize}
These three sets of rules give rise to the definition of three different spans: the strong span, direct span, and Cartesian span. Each of these spans can have a vertex version or an edge version, depending on whether Alice and Bob are required to visit all the vertices or all the edges of the graph. For edge spans, we use the same notation as in \cite{banic}:
\begin{itemize}
	\item $\sigma^{\boxtimes}_{E}(G)$ - strong edge span, corresponding to the traditional movement rules;
	\item $\sigma^{\times}_{E}(G)$ - direct edge span, corresponding to the active movement rules;
	\item $\sigma^{\square}_{E}(G)$ - Cartesian edge span, corresponding to the lazy movement rules.
\end{itemize}

In our previous work \cite{nas}, we introduced functions that represent walks in a graph, that encompass visits to all the vertices. We refer to these functions as $l$-tracks and lazy $l$-tracks, depending on the movement rules they adhered to. Here, we revisit these definitions to establish analogous concepts, specifically functions that represent walks in a graph that visit all the edges.

We will use standard definitions and concepts of graph theory, which can be found, for example, in \cite{gross}. By $\mathbb N_l$ we denote the set of first $l$ natural numbers, that is, the set $\{1, \ldots, l \}$. For the purpose of our study, we restricted our focus to finite connected graphs. Hence, throughout the remainder of the paper, the term 'graph' specifically refers to a finite connected graph.

\begin{definition}
	Let $G=(V,E)$ be a graph, and $l\in\mathbb{N}$. We say that a surjective function $f_{l}:\mathbb{N}_{l}\longrightarrow V(G)$ is an\textbf{ $l$-track} on $G$ if $f(i)f(i+1)\in E(G)$ holds, for each $i\in\mathbb{N}_{l-1}$.
\end{definition}

\begin{definition}\label{sweep}
	Let $G=(V,E)$ be a graph, and $l\in\mathbb{N}$. We say that a surjective function $f_{l}:\mathbb{N}_{l}\longrightarrow V(G)$ is a\textbf{ lazy $l$-track} on $G$ if $f(i)f(i+1)\in E(G)$ or $f(i)=f(i+1)$ holds, for each $i\in\mathbb{N}_{l-1}$.
\end{definition}

Now we define the functions that represent walks in a graph that visit all the edges.

\begin{definition}
	Let $G=(V,E)$ be a graph, and $l\in\mathbb{N}$. We say that an $l$-track $f$ is an\textbf{ $l$-sweep} on $G$ if for each $uv\in E(G)$ there exists $i\in\mathbb{N}_{l-1}$ such that $f(i)f(i+1)=uv$.
\end{definition}

\begin{definition}\label{lazy}
	Let $G=(V,E)$ be a graph, and $l\in\mathbb{N}$. We say that a lazy $l$-track $f$ on $G$ is a \textbf{lazy $l$-sweep} on $G$ if for each $uv\in E(G)$ there exists $i\in\mathbb{N}_{l-1}$ such that $f(i)f(i+i)=uv$.
\end{definition}

 Note that an $l$-sweep is also a lazy $l$-sweep. If at some point we consider a lazy $l$-sweep, or an $l$-sweep $f$ for which the value of $l$ is unimportant, we will refer to $f$ as a lazy sweep or a sweep, to make the text more readable.
 
An analogous result to Proposition 2.3 in \cite{nas} holds, and we present it here without providing a proof.

\begin{proposition}\label{uvodna}
	Let $G$ be a graph. There exists an $l\in\mathbb{N}$, such that there exists an $l$-sweep on $G$, and for each $l'\in\mathbb{N}$, $l'>l$, there exists a lazy $l'$-sweep on $G$.
\end{proposition}

The distance between two lazy $l$-sweeps on a graph is defined analogously as the distance between two lazy $l$-tracks.

\begin{definition}
	Let $G$ be a graph and $f,g$ lazy $l$-sweeps on $G$. The \textbf{distance between $f$ and $g$} is defined as
	$$m_{G}(f,g)=\min\{d_{G}(f(i),g(i)):i\in\mathbb{N}_{l}\}.$$
\end{definition}
If it is clear to which graph $G$ we are referring to, we will write only $m(f,g)$.

\section{Edge spans of a graph}\label{edge}

We will now proceed to define the edge spans of a graph. In \cite{nas}, we introduced three variations of vertex spans, corresponding to the various movement rules that two players can adopt within a graph. To define the edge spans, we modify these definitions by replacing $l$-tracks with $l$-sweeps.\\

\begin{definition}
	Let $G$ be a graph, and $l\in\mathbb{N}$. If there is at least one lazy $l$-sweep on $G$ we say that $G$ is an $l$-\textbf{sweepable graph}.
\end{definition}

\textbf{Strong edge span}\\
Let $G$ be an $l$-sweepable graph. We define

$$M_{l}^{\boxtimes}=\max\{m_{G}(f,g):f \text{ and } g \text{ are } \text{lazy }l  \text{-sweeps on } G\}. $$

Let $G$ be a graph, and let $S\subseteq\mathbb{N}$ be the set of all integers $l$ for which $G$ is an $l$-sweepable graph. $S$ is non-empty due to Proposition \ref{uvodna}. We define

$$M^{\boxtimes}(G)=\max\{M_{l}^{\boxtimes}:l\in S\}.$$

This number is the maximal safety distance that can be kept while two players visit all the edges of a graph while following the traditional movement rules.

Analogously to \cite{nas} and \cite{banic}, it can be proven that $M^{\boxtimes}(G)$ equals $\sigma^{\boxtimes}_{E}(G)$ and from now on we use the denotation $\sigma^{\boxtimes}_{E}(G)$.\\

\textbf{Direct edge span}\\     
Let $G$ be an $l$-sweepable graph. We define

$$M^{\times}_{l}=\max\{m_{G}(f,g):f \text{ and } g \text{ are } l  \text{-sweeps on } G\}. $$

Let $G$ be a graph, and let $D\subseteq\mathbb{N}$ be the set of all integers $l$ for which $G$ is an $l$-sweepable graph. We define

$$M^{\times}(G)=\max\{M^{\times}_{l}:l\in D\}.$$

This number is the maximal safety distance that can be kept while two players visit all the edges of a graph with respect to the active movement rules. \\ 
We have $M^{\times}(G)=\sigma^{\times}_{E}(G)$.\\

\textbf{Cartesian edge span}\\      
For the Cartesian edge span, analogously as in \cite{nas}, we introduce the notion of opposite lazy $l$-sweeps.

\begin{definition}
	Let $G$ be an $l$-sweepable graph and $f,g$ lazy $l$-sweeps on $G$. We say that $f$ and $g$ are\textbf{ opposite lazy $l$-sweeps} on $G$ if
	$$f(i)f(i+1)\in E(G) \iff g(i)=g(i+1),$$
	for all $i\in\mathbb{N}_{l-1}$.
\end{definition}

Let $G$ be a graph and $l\in\mathbb{N}$ such that at least one pair of opposite lazy $l$-sweeps exists on $G$. We define
$$M^{\square}_{l}=\max\{m_{G}(f,g):f \text{ and } g \text{ are opposite lazy } l  \text{-sweeps on } G\}.$$

Let $G$ be a graph, and let $C\subseteq\mathbb{N}$ be the set of all integers $l$ for which opposite lazy $l$-sweeps exist on $G$. We define
$$M^{\square}(G)=\max\{M^{\square}_{l}:l\in C\}.$$

This number is the maximal safety distance that can be kept while two players visit all the edges of a graph with respect to the lazy movement rules.

We have $M^{\square}(G)=\sigma^{\square}_{E}(G)$.

\begin{remark}
	Note that, by Observation 3.2. from \cite{banic}, for any connected graph $G$ it follows:
	$$\sigma^{\boxtimes}_{E}(G)\leq\sigma^{\boxtimes}_{V}(G)\leq \rad G;$$
	$$\sigma^{\times}_{E}(G)\leq\sigma^{\times}_{V}(G)\leq \rad G;$$
	$$\sigma^{\square}_{E}(G)\leq\sigma^{\square}_{V}(G)\leq \rad G,$$
    where $\rad G$ is the radius of graph $G$. 
	
\end{remark}

\section{Preliminaries}\label{prelim}

When studying direct spans in a graph, we will frequently construct $l$-tracks or $l$-sweeps, $l \in \mathbb{N}$, and when examining Cartesian spans, it will be necessary to modify these functions to become opposite lazy $k$-tracks or opposite lazy $k$-sweeps, $k \in \mathbb{N}$, while still adhering to the same sequence of vertices.
The concept of $l$-tracks or $l$-sweeps can be visualized as Alice and Bob following active movement rules within a graph, visiting all the vertices or all the edges. However, if they are instead required to abide by lazy movement rules, these $l$-tracks or $l$-sweeps will be adjusted so that Alice and Bob visit the same sequence of vertices in their walks. In each step, only one of them moves, and they alternate their movements.
In fact, we have already used this in Theorem 4.2. in \cite{nas}. We defined two opposite lazy $(2l-1)$-tracks using two $l$-tracks. Here we see the need to define those kinds of functions formally, to facilitate smoother and more coherent proofs in the upcoming sections. Therefore, we now introduce the notion of opposite lazy $(2l-1)$-tracks induced by $l$-tracks, and opposite lazy $(2l-1)$-sweeps induced by $l$-sweeps.

\begin{definition}
	Let $G$ be a graph, and $f,g$ be $l$-tracks on $G$. We say that $f',g':\mathbb{N}_{2l-1}\to V(G)$ defined with
	$$f'(i)=f\left(\left\lceil\frac{i+1}{2}\right\rceil\right);$$
	$$g'(i)=g\left(\left\lceil\frac{i}{2}\right\rceil\right).$$
	are \textbf{opposite lazy $(2l-1)$-tracks induced by $l$-tracks $f$ and $g$.}
\end{definition}
\begin{remark}
	Note that in Theorem 4.2. in \cite{nas}, we prove that $f'$ and $g'$ defined in this way are indeed opposite lazy $(2l-1)$-tracks and so here we just define them formally. Moreover, from the same proof it follows that $m(f',g')\geq m(f,g)-1$.
\end{remark}
Analogously, we define \textbf{opposite lazy $(2l-1)$-sweeps induced by $l$-sweeps}.

In the following two sections, we will utilize the concept of the shortest walk in a complete graph, denoted by $K_n$ where $n\geq 2$, and in a complete bipartite graph denoted as $K_{2,n-2}$ where $n\geq 4$, which traverses all the edges of the respective graphs. Therefore, we present the preliminary results here.
It is important to note that the length of a walk in a graph is determined by the number of edges in that walk. For odd values of $n$, the result for the complete graph $K_n$ is trivial since the graph is Eulerian, while the result for even $n$ is based on observations from the Chinese Postman problem and the algorithms used to solve it. The proof that the given number represents the length of the shortest walk can be found in algorithms and claims presented in \cite{kinez}. Here, we provide a brief explanation of how the walk is obtained.

\begin{proposition}\label{Euler}
	Let $n\geq 2$. The length $s_{K_n}$ of the shortest walk that visits every edge at least once in a graph $K_n$ is equal to:
	$$s_{K_n}=\left\{ \begin{array}{lc}
		\frac{n^2-n}{2}, & n \text{ odd;} \\
		\frac{n^2-2}{2}, & n \text{ even.}
	\end{array} \right.$$
\end{proposition}
\begin{proof}
	For odd $n$, every $K_n$ has an Eulerian circuit, a trail that visits all the edges exactly once and starts and ends in the same vertex. The length of such a walk is $\frac{n^2-n}{2}$, since that is the size of the graph $K_n$.\\
	For $n=2$, $K_2$ has an Eulerian trail, and its length is $1$.\\
	For even $n$, $n\geq 4$, $K_n$ does not have an Eulerian trail since it has more than $2$ vertices of odd degree, so we ask ourselves what is the minimal number of edges a walk must visit more than once in order to traverse all the edges. Note that we do not ask for such a walk to be closed. Consider a partition of vertices of $K_n$ into pairs, and let $G$ be a multigraph obtained from $K_n$, by adding one edge between two vertices in each but one pair in the observed partition. There are only $2$ vertices in $G$ that have an odd degree. It is easy to see that the number of edges added is $\frac{n-2}{2}$, so $E(G)=\binom{n}{2}+\frac{n-2}{2}=\frac{n^2-2}{2}$. $G$ has an Eulerian trail of length $\frac{n^2-2}{2}$ which corresponds to a walk in $K_n$ in which $\frac{n-2}{2}$ edges are traversed twice, so it follows $s_{K_n}=\frac{n^2-2}{2}$, for even $n\geq 4$. 
\end{proof}

\begin{corollary}\label{duljina}
	For each $n\in\mathbb{N}$, there exists a $k$-sweep on $K_n$, where 
	$$k=\left\{ \begin{array}{lc}
		\frac{n^2-n+2}{2}, & n \text{ odd;} \\
		\frac{n^2}{2}, & n \text{ even.}
	\end{array} \right.$$
\end{corollary}
\begin{proof} For $n=1$ the result is trivial and for $n\geq 2$ the claim follows directly from the shortest walks in Proposition \ref{Euler}.

\end{proof}

The result for a complete bipartite graph $K_{2,n-2}$, $n\geq 4$, is given in the subsequent claim.

\begin{proposition}\label{bipartitni}
	There exists a $(2n-3)$-sweep on $K_{2,n-2}$, for each $n\geq 4$.
\end{proposition}
\begin{proof}
	It is easy to see that $K_{2,n-2}$ has $2n-4$ edges.
  For even $n$, all the vertices have an even degree and for odd $n$ there are exactly two vertices of an odd degree, so in each case the graph has an Eulerian trail and the length of the shortest walk that visits every edge at least once in $K_{2,n-2}$, is $2n-4$. We conclude that a $(2n-3)$-sweep exists.

\end{proof}

\section{Results on edge spans}\label{res_edge}

Let us recap the results for the edge spans given in \cite{banic}. 
\begin{itemize}
	\item If $\rad G=1$ then $\sigma^{\boxtimes}_{V}(G)=\sigma^{\boxtimes}_{E}(G)$, $\sigma^{\times}_{V}(G)=\sigma^{\times}_{E}(G)$ and $\sigma^{\square}_{V}(G)=\sigma^{\square}_{E}(G)$;
	\item For a path $P_n$ we have $\sigma^{\boxtimes}_{V}(G)=\sigma^{\boxtimes}_{E}(G)=1$ and $\sigma^{\square}_{V}(G)=\sigma^{\square}_{E}(G)=0$;
	\item If $G$ is an $n$-friendly graph (a family of graphs defined in \cite{banic}, Definition 6.2), and $n=\rad (G)$, then $$\sigma^{\boxtimes}_{V}(G)=\sigma^{\boxtimes}_{E}(G)=\rad G,$$ 
	$$\sigma^{\times}_{V}(G)=\sigma^{\times}_{E}(G)=\rad G.$$
\end{itemize}

In paper \cite{nas}, various results were presented concerning the relationship between different vertex spans, along with specific vertex span values for various graph classes. In this study, our focus shifts to determining the values of edge spans for certain graph classes and exploring the relationship between different edge spans.

We begin with the relation between strong, direct and Cartesian edge span, analogous to Proposition 4.1. in \cite{nas}. We state the result without proof since it is analogous to the vertex variant.

\begin{proposition}\label{strong_max}
    Let $G$ be a graph. We have $\sigma^{\boxtimes}_{E}(G)\geq\max \{\sigma^{\times}_{E}(G),\sigma^{\square}_{E}(G)\}$.
\end{proposition}

Let us quickly demonstrate that the values of vertex spans and edge spans are equal for certain classes of graphs.

\begin{proposition}\label{stablo}
	If $G$ is a tree, then $\sigma^{\boxtimes}_{V}(G)=\sigma^{\boxtimes}_{E}(G)$, $\sigma^{\times}_{V}(G)=\sigma^{\times}_{E}(G)$ and $\sigma^{\square}_{V}(G)=\sigma^{\square}_{E}(G)$.
\end{proposition}
\begin{proof}
	These claims follow from the fact that, in a tree, there is a unique path between every two vertices, so each $l$-track is also an $l$-sweep, and each lazy $l$-track also a lazy $l$-sweep.
\end{proof}

\begin{proposition}\label{ciklus_edge}
	If $n\geq 3$ then $\sigma^{\boxtimes}_{V}(C_n)=\sigma^{\boxtimes}_{E}(C_n)$, $\sigma^{\times}_{V}(C_n)=\sigma^{\times}_{E}(C_n)$ and $\sigma^{\square}_{V}(C_n)=\sigma^{\square}_{E}(C_n)$.
\end{proposition}
\begin{proof}
	All claims follow from the fact that a cycle $C_n$ is isomorphic to its line graph $L(C_n)$. Obviously, $\sigma^{*}_{E}(C_n)=\sigma^{*}_{V}(L(C_n))$, $*\in\{\boxtimes,\times,\square\}$. (The facts that $\sigma^{\boxtimes}_{V}(C_n)=\sigma^{\boxtimes}_{E}(C_n)$ and $\sigma^{\times}_{V}(C_n)=\sigma^{\times}_{E}(C_n)$ also follow from the results in \cite{banic}, because $C_n$ is a $k$-friendly graph for $n\geq 2k$.)
\end{proof}

In all of the aforementioned families, as well as the previously studied ones, the strong edge span is consistently equal to the larger value between the direct edge span and the Cartesian edge span. However, we will now introduce a family of graphs where this equality does not hold.\\
Let $G=(V,E)$ be a graph. An \textbf{edge subdivision} of the edge $uv\in E(G)$ is the operation of replacing $uv$ with a path $uwv$, through a new vertex $w$. The result of this operation is a graph $H$:
$$H=(V(G)\cup \{w\},(E(G)\setminus\{uv\})\cup \{uw,wv\}).$$
A graph which has been derived from $G$ by a series of edge subdivisions is called a \textbf{$G$-subdivision} or a \textbf{subdivision of} $G$.

Let $n\geq 4$ and let $K_n^+$ be a $K_n$-subdivision obtained by exactly one edge subdivision. Graphs $K_5^+$ and $K_6^+$ are shown in Figure \ref{Kn+}. 

\begin{center}
	\begin{tikzpicture}
		\coordinate (v1) at (306:2cm);
		\coordinate (v2) at (18:2cm);
		\coordinate (v3) at (90:2cm);
		\coordinate (v4) at (162:2cm);
		\coordinate (v5) at (234:2cm);
		\coordinate (v6) at ($(v1)!0.5!(v5)$);
		\node (v1) at (v1) [point] {};
		\node [below=1pt]  at (v1) {$v_1$};
		\node (v2) at (v2) [point] {};
		\node [right=1pt]  at (v2) {$v_2$};
		\node (v3) at (v3) [point] {};
		\node [above=1pt]  at (v3) {$v_3$};
		\node (v4) at (v4) [point] {};
		\node [left=1pt]  at (v4) {$v_4$};
		\node (v5) at (v5) [point] {};
		\node [below=1pt]  at (v5) {$v_5$};
		\node (v6) at (v6) [point] {};
		\node [below=1pt]  at (v6) {$w$};
		\draw (v1) -- (v2) -- (v3) -- (v4) -- (v5) -- (v6) -- (v1) -- (v3) -- (v5) -- (v2) -- (v4) -- (v1);
	\end{tikzpicture}
	\hspace{20pt}
	\begin{tikzpicture}
		\coordinate (v1) at (300:2cm);
		\coordinate (v2) at (0:2cm);
		\coordinate (v3) at (60:2cm);
		\coordinate (v4) at (120:2cm);
		\coordinate (v5) at (180:2cm);
		\coordinate (v6) at (240:2cm);
		\coordinate (v7) at ($(v1)!0.5!(v6)$);
		\node (v1) at (v1) [point] {};
		\node [below=1pt]  at (v1) {$v_1$};
		\node (v2) at (v2) [point] {};
		\node [right=1pt]  at (v2) {$v_2$};
		\node (v3) at (v3) [point] {};
		\node [above=1pt]  at (v3) {$v_3$};
		\node (v4) at (v4) [point] {};
		\node [above=1pt]  at (v4) {$v_4$};
		\node (v5) at (v5) [point] {};
		\node [left=1pt]  at (v5) {$v_5$};
		\node (v6) at (v6) [point] {};
		\node [below=1pt]  at (v6) {$v_6$};
		\node (v7) at (v7) [point] {};
		\node [below=1pt]  at (v7) {$w$};
		\draw (v1) -- (v2) -- (v3) -- (v4) -- (v5) -- (v6) -- (v7) -- (v1) (v1) -- (v3) -- (v5) -- (v1) (v4) -- (v6) -- (v2) -- (v4) (v1) -- (v4) (v2) -- (v5) (v3) -- (v6);
	\end{tikzpicture}
	\captionof{figure}{$K_5^+$ and $K_6^+$} 
	\label{Kn+}
\end{center}

The next theorem shows that for the graph $K_n^+$, $n\geq 4$, we have
$$\sigma^{\boxtimes}_{E}(K_n^+)>\max \{\sigma^{\times}_{E}(K_n^+),\sigma^{\square}_{E}(K_n^+)\}.$$
This result is intriguing because, currently, no graph is known to exist where an analogous claim holds for vertex spans. Specifically, there is no known graph for which the strong vertex span is not equal to the maximum value between the direct and Cartesian vertex spans, for all considered graphs thus far. Furthermore, it is worth noting that the value of the direct edge span for $K_n^+$ does not equal the value of its direct vertex span. We will subsequently establish, as demonstrated in Proposition \ref{najmanji}, that $K_4^+$ is the graph of the smallest order for which this property holds.\\

\begin{theorem}\label{tri}
	If $n\geq 4$ then
	\begin{enumerate}[label=(\roman*)]
		\item $\sigma^{\boxtimes}_{V}(K_n^+)=\sigma^{\boxtimes}_{E}(K_n^+)=2;$
		\item $\sigma^{\times}_{V}(K_n^+)=2,\quad\sigma^{\times}_{E}(K_n^+)=1;$
		\item $\sigma^{\square}_{V}(K_n^+)=\sigma^{\square}_{E}(K_n^+)=1.$
	\end{enumerate}
\end{theorem}

\begin{proof}
	Let $n\geq 4$ and denote the vertices of $K_n^+$ by $v_1,v_2,...,v_n,w$, in such a way that $w\notin V(K_n)$ and $v_1w,v_nw\in E(K_n^+)$. Note that $\rad (K_n^+)=2$ and that all the vertices in $\{v_2,...v_{n-1}\}$, are at distance $2$ from $w$ in $K_n^+$. Notice that the only vertex at distance $2$ from $v_1$ is $v_n$ and vice versa. Those are the only distances $2$ between the vertices of $K_n^+$. We can see that $$K_n^+-(\{w\}\cup N(w))\cong K_{n-2},$$ 
	with the set of vertices $\{v_2,...v_{n-1}\}$. Let us denote $H_1=K_n^+-(\{w\}\cup N(w))$. Also observe that $$K_n^+-\{w\}-E(H_1)\cong K_{2,n-2},$$ and denote $H_2=K_n^+-\{w\}-E(H_1)$. It is easy to see that $\{E(H_1),E(H_2),\{v_1w,v_nw\}\}$ is a partition of the set of edges of $K_n^+$, (Figure \ref{part}). 
	\begin{center}
		\begin{tikzpicture}
			\coordinate (v1) at (300:2cm);
			\coordinate (v2) at (0:2cm);
			\coordinate (v3) at (60:2cm);
			\coordinate (v4) at (120:2cm);
			\coordinate (v5) at (180:2cm);
			\coordinate (v6) at (240:2cm);
			\coordinate (v7) at ($(v1)!0.5!(v6)$);
			\node (v1) at (v1) [point] {};
			\node [below=1pt]  at (v1) {$v_1$};
			\node (v2) at (v2) [point] {};
			\node [right=1pt]  at (v2) {$v_2$};
			\node (v3) at (v3) [point] {};
			\node [above=1pt]  at (v3) {$v_3$};
			\node (v4) at (v4) [point] {};
			\node [above=1pt]  at (v4) {$v_4$};
			\node (v5) at (v5) [point] {};
			\node [left=1pt]  at (v5) {$v_5$};
			\node (v6) at (v6) [point] {};
			\node [below=1pt]  at (v6) {$v_6$};
			\node (v7) at (v7) [point] {};
			\node [below=1pt]  at (v7) {$w$};
			\draw [line width=2pt, blue!] (v6) -- (v7) -- (v1);
			\draw [line width=2pt, green!] (v1) -- (v2) (v1) -- (v3) (v1) -- (v4) (v1) -- (v5) (v6) -- (v2) (v6) -- (v3) (v6) -- (v4) (v6) -- (v5);
			\draw [line width=2pt, red!] (v2) -- (v3) -- (v4) -- (v5) -- (v2) (v2) -- (v4) (v3) -- (v5);
		\end{tikzpicture}
		\captionof{figure}{The partition of edges of $K_6^+$} 
		\label{part}
	\end{center}
	
	(i) Let us prove that $\sigma^{\boxtimes}_{V}(K_n^+)=\sigma^{\boxtimes}_{E}(K_n^+)=2$. It is enough to prove $\sigma^{\boxtimes}_{E}(K_n^+)\geq 2$ because $\sigma^{\boxtimes}_{E}(K_n^+)\leq\sigma^{\boxtimes}_{V}(K_n^+)\leq\rad (K_n^+)=2$. We will show that there exist lazy $l$-sweeps $f,g$ in $K_n^+$, such that $m(f,g)=2$, for some $l\in\mathbb{N}$.\\       
	
	We will construct functions $f$ and $g$ in four stages. Before delving into the formal construction, let us provide a brief description in terms of Alice and Bob, ensuring that they maintain a constant distance of $2$ from each other throughout their walks. We give an example of their movement in $K_5^+$ in Example \ref{AB}.
	In the first stage, Alice will traverse all edges of $H_1$, following the shortest walk that covers each of these edges (Proposition \ref{Euler} and Corollary \ref{duljina}), while Bob will remain stationed at the vertex $w$. 
	In the second stage, Alice will visit all the edges of $H_2$, specifically the edges of the form $v_1v_j$ and $v_{n}v_j$, where $j\in\{2,...,n-1\}$. Here is how it will be accomplished: Starting from the last vertex Alice visited in $H_1$, she moves to vertex $v_1$. Subsequently, she follows the shortest walk in $H_2$ that covers all of its edges (Proposition \ref{bipartitni}). During this stage, Bob will alternate between the vertices $w$, $v_1$, and $v_n$ according to the following pattern: Whenever Alice is at vertex $v_1$, Bob will be at vertex $v_n$. Whenever Alice is at vertex $v_j$, $j\in\{2,...,n-1\}$, Bob will be at vertex $w$. Finally, whenever Alice is at vertex $v_n$, Bob will be at vertex $v_1$. At this point, Alice will have visited all the edges of $K_n^+$ except $wv_1$ and $wv_n$, while Bob will have only visited those two edges. To ensure that both Alice and Bob visit all the edges, Bob will need to repeat Alice's steps from the first two stages, while Alice repeats Bob's steps. Consequently, in the third stage, Bob will visit all the edges of $H_2$, while Alice alternates between the vertices $w$, $v_1$, and $v_n$, thus covering the remaining edges she had not visited in the previous stages. Finally, in the last stage, Bob will visit the edges of $H_1$, while Alice will remain at the vertex $w$. In the formal construction of functions $f$ and $g$ that represent this movement, we will also determine the number of steps involved in each stage, so we have the number $l\in\mathbb{N}$, for which $f$ and $g$ are lazy $l$-sweeps in $K_n^+$ and it is 
	$$l=2k+4n-7,$$ where
	$$k=\left\{ \begin{array}{lc}
		\frac{n^2-5n+8}{2}, & n \text{ odd} \\
		\frac{n^2-4n+4}{2}, & n \text{ even.}
	\end{array} \right. $$
	Let us proceed with the formal construction.\\
	
	First stage:\\
	By Proposition \ref{Euler} and Corollary \ref{duljina}, let $W_1=(v_2,...,v_x)$, $x\in\{2,...,n-1\}$ be the shortest walk that traverses all the edges of $H_1$ and let $f$ and $g$ be functions such that 
	$$(f(1),f(2),...,f(k))=W_1 \text{ and } g(i)=w, \text{ for all } i\in\{1,...,k\},$$
	where $$k=\left\{ \begin{array}{lc}
		\frac{n^2-5n+8}{2}, & n \text{ odd} \\
		\frac{n^2-4n+4}{2}, & n \text{ even.}
	\end{array} \right. $$
	Second stage:\\
	By Proposition \ref{bipartitni}, let $W_2=(v_1,...v_y)$, where $y=1$, if $n$ is even, and $y=n$, if $n$ is odd, be the shortest walk that traverses all edges of $H_2$. Let
	$$(f(k+1),f(k+2),...,f(k+2n-3))=W_2,$$
	and for each $i\in\{k+1,...,k+2n-3\}$ let
	$$g(i)=\left\{ \begin{array}{lc}
		v_n, & f(i)=v_1; \\
		v_1, & f(i)=v_n; \\
		w, & f(i)\in V(H_1).
	\end{array} \right. $$
	
	Third stage:\\
	Observe that the following is true:
	\begin{itemize}
		\item For even $n$, it follows $f(k+2n-3)=f(k+1)=v_1$ and  $g(k+2n-3)=g(k+1)=v_n$;
		\item For odd $n$, it follows $f(k+2n-3)=g(k+1)=v_n$ and  $g(k+2n-3)=f(k+1)=v_1$. 
	\end{itemize}
	
	Since $v_1$ and $v_n$ are both adjacent to all vertices in $V(H_1)$, in both cases we can proceed with observing $W'_2=(g(k+2n-3),...,v_y)$, $y\in\{1,n\}$, the shortest walk that traverses all the edges of $H_2$. We define
	$$(g(k+2n-3),...,g(k+4n-7))=W'_2,$$
	and for each $i\in\{k+2n-2,...,k+4n-7\}$ let
	$$f(i)=\left\{ \begin{array}{lc}
		v_n, & g(i)=v_1; \\
		v_1, & g(i)=v_n; \\
		w, & g(i)\in V(H_1).
	\end{array} \right. $$
	Fourth stage:\\
	For the last part, let
	$$(g(k+4n-6),g(k+4n-5),...,g(2k+4n-7))=W_1,$$  
	and
	$$f(i)=w, \text{ for all } i\in\{k+4n-6,...,2k+4n-7\}.$$
	
	In each stage of the construction, it is clear that $m(f,g)=2$, and with the observed number of steps we conclude that $f$ and $g$ are lazy $(2k+4n-7)$-sweeps in $K_n^+$, which proves $\sigma^{\boxtimes}_{V}(K_n^+)=\sigma^{\boxtimes}_{E}(K_n^+)=2$.\\
	
	
	(ii) Let us show that $\sigma^{\times}_{V}(K_n^+)=2$ and $\sigma^{\times}_{E}(K_n^+)=1$.
	From the construction of lazy $(2k+4n-7)$-sweeps $f$ and $g$ in (i), it is clear that their restrictions to the set $\{k+1,...,k+4n-7\}$ (the second and third stage), are $(4n-7)$-tracks in $K_n^+$ that keep the distance $2$, so $\sigma^{\times}_{V}(K_n^+)=2$.\\
	Let $f$ and $g$ be $l$-sweeps in $K_n^+$, for any $l\in\mathbb{N}$. We will demonstrate that $m(f,g)<2$. Let $v_iv_j$ be any edge in $E(H_1)$, $i,j\in\{2,3,...,n-1\}$. Since $f$ is an $l$-sweep, there exists some $k<l$ such that $v_iv_j=f(k)f(k+1)$. Since $w$ is the only vertex at the distance $2$ from both $v_i$ and $v_j$, we may assume $g(k)=w$ or $g(k+1)=w$, however, $g(k)\neq g(k+1)$, so we conclude $m(f,g)<2$. This proves $\sigma^{\times}_{E}(K_n^+)<2$ and since it is easy to see that $\sigma^{\times}_{E}(K_n^+)>0$, for $n\geq 4$, we conclude $\sigma^{\times}_{E}(K_n^+)=1$.\\


	(iii)  Let us show that $\sigma^{\square}_{V}(K_n^+)=\sigma^{\square}_{E}(K_n^+)=1$. It is easy to see that $\sigma^{\square}_{V}(K_n^+),\sigma^{\square}_{E}(K_n^+)\geq 1$ and since $\sigma^{\square}_{V}(K_n^+)\geq\sigma^{\square}_{E}(K_n^+)$, it remains to show that $\sigma^{\square}_{V}(K_n^+)<2$. All neighbors of $v_1$ and $v_n$ are the same, i.e. $N(v_1)=N(v_n)$, and $v_n$ is the only vertex on the distance $2$ from $v_1$. So for any opposite lazy $l$-tracks $f,g$, for which $m(f,g)=2$, if $f(k)=v_1$, for some $k<l$, then $g(k)=v_n$ must hold. Assume $f(k+1)=f(k)$. Now $g(k+1)$ is a vertex that is on the distance $1$ from $v_1$, so $m(f,g)<2$, which proves the claim.
\end{proof}

\begin{example}
	To make the construction in (i) of Theorem \ref{tri} clearer, let us give an example of Alice and Bob's movement in $K_5^+$ in Table \ref{AB}. Graph $K_5^+$ is shown in Figure \ref{Kn+}. We have $k=4$ and $2k+4n-7=21$. The third stage includes steps 12-17, however, to emphasize the symmetry between Alice and Bob's movements in the first two and the last two stages, we repeat the 11th step in the third row of the table.
	\captionof{table}{Four stages of Alice and Bob's movement in $K_5^+$}\label{AB} 
	\begin{tabular}{ |c|p{0.4cm}|p{0.4cm}|p{0.4cm}|p{0.4cm}| } 
		\hline
		& 1 & 2 & 3 & 4 \\ 
		\hline
		Alice & $v_{2}$ & $v_{3}$ & $v_{4}$ & $v_{2}$  \\ 
		Bob & $w$ & $w$ & $w$ & $w$  \\ 
		\hline
	\end{tabular}\\
	\begin{tabular}{ |c|p{0.4cm}|p{0.4cm}|p{0.4cm}|p{0.4cm}|p{0.4cm}|p{0.4cm}|p{0.4cm}| } 
		\hline
		& 5 & 6 & 7 & 8 & 9 & 10 & 11\\ 
		\hline
		Alice & $v_{1}$ & $v_{4}$ & $v_{5}$ & $v_{3}$ & $v_{1}$ & $v_{2}$ & $v_5$ \\ 
		Bob & $v_{5}$ & $w$ & $v_{1}$ & $w$ & $v_{5}$ & $w$ & $v_{1}$ \\ 
		\hline
	\end{tabular}\\
	\begin{tabular}{ |c|p{0.4cm}|p{0.4cm}|p{0.4cm}|p{0.4cm}|p{0.4cm}|p{0.4cm}|p{0.4cm}| } 
		\hline
		& 11 & 12 & 13 & 14 & 15 & 16 & 17 \\ 
		\hline
		Alice & $v_5$ & $w$ & $v_{1}$ & $w$ & $v_{5}$ & $w$ & $v_{1}$ \\ 
		Bob & $v_1$ & $v_{4}$ & $v_{5}$ & $v_{3}$ & $v_{1}$ & $v_{2}$ & $v_5$ \\ 
		\hline
	\end{tabular}\\
	\begin{tabular}{ |c|p{0.4cm}|p{0.4cm}|p{0.4cm}|p{0.4cm}| } 
		\hline
		& 18 & 19 & 20 & 21 \\ 
		\hline
		Alice & $w$ & $w$ & $w$ & $w$  \\ 
		Bob & $v_{2}$ & $v_{3}$ & $v_{4}$ & $v_{2}$   \\ 
		\hline
		
	\end{tabular}

\end{example}

Theorem \ref{tri} proves the existence of a graph $G$ for which $$\sigma^{\boxtimes}_{E}(G)>\max \{\sigma^{\times}_{E}(G),\sigma^{\square}_{E}(G)\}.$$
It remains an open problem whether there exists a graph $H$ for which $\sigma^{\boxtimes}_{V}(H)>\max \{\sigma^{\times}_{V}(H),\sigma^{\square}_{V}(H)\}$, or the value of strong vertex span always equals the maximal value of the direct and Cartesian vertex span. Before we prove a general relation between strong, direct and Cartesian edge span, we show the following.

\begin{proposition}\label{najmanji}
	Graph $G=K_4^+$ is the graph of minimal order and size for which $\sigma^{\times}_{V}(G)\neq\sigma^{\times}_{E}(G)$.
\end{proposition}
\begin{proof}
	Every connected graph of order $3$ or $4$ is a cycle or a tree or has the radius $1$ so, by the previous results, the value of its direct edge span is equal to the value of its direct vertex span.
	Connected graphs of order $5$, with size $4$, are trees, so it is left to observe all connected graphs $G$ of order $5$, with size $5$ or $6$, such that $G\neq C_5$ and $\rad (G)>1$, which means $\Delta (G)=3$. All such graphs are shown in Figure \ref{prim} with their direct span values given.

	\begin{center}
		\begin{tabular}{ccc}
			\begin{tikzpicture}[scale=0.8]
				\node (1) at (0,0) [point] {};
				\node (2) at (2,0) [point] {};
				\node (3) at (2,2) [point] {};
				\node (4) at (0,2) [point] {};
				\node (5) at (1,3) [point] {};
				\draw (1) -- (4) -- (5) -- (3) -- (2) (4) -- (3); 
			\end{tikzpicture} 
			&
			\begin{tikzpicture}[scale=0.8]
				\node (1) at (0,0) [point] {};
				\node (2) at (2,0) [point] {};
				\node (3) at (2,2) [point] {};
				\node (4) at (0,2) [point] {};
				\node (5) at (1,3) [point] {};
				\draw (2) -- (1) -- (4) -- (5) -- (3) -- (4); 
			\end{tikzpicture}
			& 
			\begin{tikzpicture}[scale=0.8]
				\node (1) at (0,0) [point] {};
				\node (2) at (2,0) [point] {};
				\node (3) at (2,2) [point] {};
				\node (4) at (0,2) [point] {};
				\node (5) at (1,3) [point] {};
				\draw (5) -- (4) -- (1) -- (2) -- (3) -- (4);
			\end{tikzpicture}\\
			$\sigma^{\times}_V(G) = 1$ &
			$\sigma^{\times}_V(G) = 1$ &
			$\sigma^{\times}_V(G) = 2$ \\
			$\sigma^{\times}_E(G) = 1$ &
			$\sigma^{\times}_E(G) = 1$ &
			$\sigma^{\times}_E(G) = 2$ \\
			\begin{tikzpicture}[scale=0.8]
				\node (1) at (0,0) [point] {};
				\node (2) at (2,0) [point] {};
				\node (3) at (2,2) [point] {};
				\node (4) at (0,2) [point] {};
				\node (5) at (1,3) [point] {};
				\draw (4) -- (1) -- (2) -- (3) -- (4) -- (5) -- (3);
			\end{tikzpicture}
			&
			\begin{tikzpicture}[scale=0.8]
				\node (1) at (0,0) [point] {};
				\node (2) at (2,0) [point] {};
				\node (3) at (2,2) [point] {};
				\node (4) at (0,2) [point] {};
				\node (5) at (1,3) [point] {};
				\draw (5) -- (4) -- (3) -- (2) -- (1) -- (4) (1) -- (3);
			\end{tikzpicture}
			&
			\begin{tikzpicture}[scale=0.8]
				\node (1) at (0,0) [point] {};
				\node (2) at (2,0) [point] {};
				\node (3) at (2,2) [point] {};
				\node (4) at (0,2) [point] {};
				\node (5) at (1,1) [point] {};
				\draw (1) -- (2) -- (3) -- (4) -- (1) -- (5) -- (3);
			\end{tikzpicture}\\
			$\sigma^{\times}_V(G) = 2$ &
			$\sigma^{\times}_V(G) = 1$ &
			$\sigma^{\times}_V(G) = 2$\\
			$\sigma^{\times}_E(G) = 2$ &
			$\sigma^{\times}_E(G) = 1$ &
			$\sigma^{\times}_E(G) = 2$
		\end{tabular}
		\captionof{figure}{Connected graphs of order $5$ with size $5$ or $6$ and $\Delta (G)=3$} 
		\label{prim}
	\end{center}
	
	Subsequently, a graph $G$ for which $\sigma^{\times}_{V}(G)\neq\sigma^{\times}_{E}(G)$ has to have at least $5$ vertices and at least $7$ edges, and $K_4^+$ is one such graph.
\end{proof}

Now, let us turn our attention to the relation between different edge spans in general. First, we prove that the direct and Cartesian edge spans of any graph differ at most $1$ (an analogous result to Theorem 4.2. in \cite{nas}).

\begin{theorem}\label{analog}
Let $G$ be a graph. Then $$|\sigma^{\times}_E(G)-\sigma^{\square}_E(G)| \leq 1.$$
\end{theorem}
\begin{proof}
    For a trivial graph, the claim is obvious. Let $G$ be a non-trivial graph. The claim will be shown in two steps:
    \begin{itemize}
    \item [(i)]$ \sigma^{\square}_E(G) \geq \sigma^{\times}_E(G)-1$;
    \item [(ii)]$ \sigma^{\times}_E(G) \geq \sigma^{\square}_E(G)-1$.
\end{itemize} 
(i) Let $\sigma^{\times}_E(G)=r$ and let $f,g$ be $l$-sweeps for which this span value is achieved. Let $f',g'$ be opposite lazy $(2l-1)$-sweeps induced by $l$-sweeps $f$ and $g$. Notice that by definition, each lazy $l$-sweep is also a lazy $l$-track. By the proof of Theorem 4.2. in \cite{nas} it now follows that
$$m_G(f',g')\geq r-1,$$ hence
$$ \sigma^{\square}_E(G) \geq \sigma^{\times}_E(G)-1.$$
(ii) Let $\sigma^{\square}_E(G)=r$ and let $f,g$ be opposite lazy $l$-sweeps for which this span value is achieved. $f$ and $g$ are by definition opposite lazy $l$-tracks, so by the proof of Theorem 4.2. in \cite{nas} we know that for some $k\in\mathbb{N}$ there exist $k$-tracks $f',g'$ such that $m_G(f',g')\geq r-1$. Since $f$ is a sweep, it is clear that for each edge $e\in E(G)$, there will exist some $i\in\mathbb{N}$ such that $f(i)f(i+1)=e$ and from the construction of $f'$ is clear that there exist some $j\in\mathbb{N}$ such that $f'(j)f'(j+1)=f(i)f(i+1)=e$, so $f'$ is a $k$-sweep, and analogously, so is $g'$, which proves the claim.
\end{proof}

We know from Proposition \ref{strong_max} that the strong edge span is greater than or equal to the largest of the values of direct and Cartesian spans, and we aim to improve that result by finding the upper bound for the strong edge span. We begin by proving that the direct edge span is greater than or equal to the value of the strong edge span minus $1$.

\begin{theorem}\label{strong_direct}
Let $G$ be a graph. Then 
     $$\sigma^{\times}_{E}(G) \geq \sigma^{\boxtimes}_{E}(G)-1.$$
\end{theorem}
\begin{proof}
    Let $\sigma^{\boxtimes}_{E}(G)=d$. It follows that there is some $l\in\mathbb{N}$ and two lazy $l$-sweeps $f$ and $g$ such that $m_G(f,g)=d$, that is, for each $i\in \{1,...,l\}, d_G(f(i),g(i))\geq d$. Note that for each pair $(i,j)$ such that $|i-j|=1$, we have $d(f(i),g(j)) \geq d-1$. We will construct two $k$-sweeps, $f'$ and $g'$, for some $k\in\mathbb{N}$, such that $m_G(f',g')\geq d-1$, which will prove the claim. First, to make the proof easier to read, we will say that the lazy sweep $f$ moves in step $i$ if $f(i)\neq f(i+1)$ and that the lazy sweep $f$ stands still in step $i$ if $f(i)= f(i+1)$. Each lazy $l$-sweep generates an $(l-1)$-sequence $(x_1,...,x_{l-1})$, $x_i\in\{M, S\}$, $i=1,...,l-1$, which we will refer to as the movement sequence, in the following way:
$$x_i =
\left\{
	\begin{array}{ll}
		M & \mbox{if lazy sweep moves in step } i; \\
		S & \mbox{if lazy sweep stands still in step } i.
	\end{array}
\right.$$
Let $X=(x_1,...,x_{l-1})$ and $Y=(y_1,...,y_{l-1})$ be movement sequences of $f$ and $g$, respectively. If there exists an $i\in \{1,...,l-1\}$ such that $x_i=y_i=S$, that is, $f(i)=f(i+1)$ and $g(i)=g(i+1)$, then $f|_{\mathbb{N}_l\setminus\{i\}},g|_{\mathbb{N}_l\setminus\{i\}}$, are also lazy $(l-1)$-sweeps that keep the minimal distance $d$, so we can assume that for each $i\in \{1,...,l-1\}$ at least one of $x_i$ and $y_i$ is equal to $M$. For our construction, consider the pair $(X,Y)$, and define a new sequence $Z=(z_1,...,z_{l-1})$, which we will refer to as the activity sequence, where for each $i\in\{1,...,l-1\}$:
$$z_i =
\left\{
	\begin{array}{ll}
		A & \mbox{if } (x_i,y_i)=(M,M); \\
		P & \mbox{if } (x_i,y_i)=(M,S) \text{ or } (x_i,y_i)=(S,M).
	\end{array}
\right.$$
Now $Z$ is a sequence generated by pair $(f,g)$, in which we distinguish the steps where both $f$ and $g$ move (denoted by $A$, as active step) and steps where exactly one of $f$ and $g$ moves (denoted by $P$, as passive step). We now continue with a partition of $Z$ into particular subsequences, which we will call blocks, and then we describe how each block generates steps for desired $k$-sweeps $f'$ and $g'$.\\
Starting from the beginning of the sequence $Z$, we inductively construct a unique partition into blocks in the following way:\\

\noindent (i) If $z_1=A$, then the subsequence $Z_1=(z_1)$ is the first block;\\
 If $z_1=P$, then:\\
\hspace*{2em} - If there exist $i\neq 1$ such that $z_i=P$ then let $i'$ be the smallest such index and let $Z_1=(z_1,...,z_{i'})$ be the first block;\\
\hspace*{2em} - If there is no $i\neq 1$ such that $z_i=P$ then $Z_1=(z_1,...,z_{l-1})$ is the only block;\\
(ii) Let $Z_1,...,Z_s$ be already defined blocks, and let $z_t$ be the last element of block $Z_s$. Then:\\
\hspace*{2em} - If $z_{t+1}=A$ then $Z_{s+1}=(z_{t+1})$;\\
\hspace*{2em} - If $z_{t+1}=P$ then:\\
\hspace*{4em} - If there exist $i>t+1$ such that $z_i=P$ then let $i'$ be the smallest such index and let $Z_{s+1}=(z_{t+1},...,z_{i'})$;\\
\hspace*{4em} - If there does not exist $i>t+1$ such that $z_i=P$ then let $Z_{s+1}=(z_{t+1},...,z_{l-1})$ be the last block.\\

This construction partitioned the sequence $Z$ into blocks such that each block is a sequence of exactly one of the following types:
$$(A)$$
$$(P,P)$$
$$(P,A,...,A,P)$$
$$(P)$$
$$(P,A,...,A)$$
Notice that blocks of type $(P)$ and $(P,A,...,A)$ may only be the last block in $Z$.
For each of these types, we will now describe what subtypes regarding the pair $(X,Y)$ they encompass, and for each subtype, we will describe the construction of steps for $f'$ and $g'$. The subtypes will be displayed as matrices, the first row showing the subsequence of $X$ and the second row for $Y$. For each block, we will use the general notation for steps, $i$ and $j$: The first element, $z_i$, of the observed block corresponds to $(x_i,y_i)$, which corresponds to steps $f(i)\rightarrow f(i+1)$ and $g(i)\rightarrow g(i+1)$ of $f$ and $g$. This generates the steps $f'(j)\rightarrow f'(j+1)$ and $g'(j)\rightarrow g'(j+1)$ of $f'$ and $g'$.
The complete $f'$ and $g'$ are obtained by consecutively applying the construction for each block of the partition, for given $f$ and $g$.

\begin{itemize}[leftmargin=*]
    \item Type $(A)$ is of the type of movement:\\
        $\begin{pmatrix}
        M \\
        M \end{pmatrix}$ which means we have
        \begin{tabular}{c}
        $f(i) \neq f(i+1)$ \\
        $g(i) \neq g(i+1)$ \\
        \end{tabular}, so we define:\\ 

        $$\equaltog{f'(j)}{f(i)} \neq \equaltog{f'(j+1)}{f(i+1)}$$
        $$\equalto{g'(j)}{g(i)} \neq \equalto{g'(j+1)}{g(i+1)}$$
    \item Type $(P,P)$ is one of the following four types of movement:
    \begin{itemize}
        \item 
        [$\begin{pmatrix}
        M & S \\
        S & M \end{pmatrix}$] which means we have
        \begin{tabular}{c}
        $f(i) \neq f(i+1) = f(i+2)$ \\
        $g(i) = g(i+1) \neq g(i+2)$ \\
        \end{tabular}, so we define:\\

        $$\equaltog{f'(j)}{f(i)} \neq \equaltog{f'(j+1)}{f(i+2)}$$ 
        $$\equalto{g'(j)}{g(i)} \neq \equalto{g'(j+1)}{g(i+2)}$$

        \item
        [$\begin{pmatrix}
        M & M\\
        S & S\end{pmatrix}$] which means we have
        \begin{tabular}{c}
        $f(i) \neq f(i+1) \neq f(i+2)$ \\
        $g(i) = g(i+1) = g(i+2) $ \\
        \end{tabular}, so we define:\\ 

        $$\equaltog{f'(j)}{f(i)} \neq \equaltog{f'(j+1)}{f(i+1)} \neq \equaltog{f'(j+2)}{f(i+2)}$$ 
        $$\equalto{g'(j)}{g(i)} \neq \equalto{g'(j+1)}{\text{any neighbor of } g(i+1)} \neq \equalto{g'(j+2)}{g(i+2)}\\$$

        \item
        [$\begin{pmatrix}
        S & M\\
        M & S\end{pmatrix}$] is analogous to type $\begin{pmatrix}
        M & S\\
        S & M\end{pmatrix}$ if we swap the roles of $f$ and $g$.

        \item
        [$\begin{pmatrix}
        S & S\\
        M & M\end{pmatrix}$] is analogous to type $\begin{pmatrix}
        M & M\\
        S & S\end{pmatrix}$ if we swap the roles of $f$ and $g$.
    \end{itemize}
    \item Type $(P,A,...,A,P)$ is one of the following four types of movement:
    \begin{itemize}
        \item
        [$\begin{pmatrix}
        M & M... & S\\
        S & M... & M\end{pmatrix}$] which means we have
        \begin{tabular}{c}
        $f(i) \neq f(i+1) \neq ... \neq f(i+k-1) = f(i+k)$ \\
        $g(i) = g(i+1) \neq ... \neq g(i+k-1) \neq g(i+k)$ \\
        \end{tabular}, so we define:\\ 

        $$\equaltog{f'(j)}{f(i)} \neq \equaltog{f'(j+1)}{f(i+1)} \neq ... \neq \equaltog{f'(j+k-2)}{f(i+k-2)} \neq \equaltog{f'(j+k-1)}{f(i+k)}$$ 
        $$\equalto{g'(j)}{g(i)} \neq \equalto{g'(j+1)}{g(i+2)} \neq ... \neq \equalto{g'(j+k-2)}{g(i+k-1)} \neq \equalto{g'(j+k-1)}{g(i+k)}$$ 

        \item
        [$\begin{pmatrix}
        M & M... & M\\
        S & M... & S\end{pmatrix}$] which means we have
        \begin{tabular}{c}
        $f(i) \neq f(i+1) \neq ... \neq f(i+k-1) \neq f(i+k)$ \\
        $g(i) = g(i+1) \neq ... \neq g(i+k-1) = g(i+k)$ \\
        \end{tabular}, so we define:\\ 

        $$\equaltog{f'(j)}{f(i)} \neq \equaltog{f'(j+1)}{f(i+1)} \neq ... \neq \equaltog{f'(j+k-2)}{f(i+k-2)} \neq \equaltog{f'(j+k-1)}{f(i+k-1)} \neq \equaltog{f'(j+k)}{f(i+k)}$$ 
        $$\equalto{g'(j)}{g(i)} \neq \equalto{g'(j+1)}{g(i+2)} \neq ... \neq \equalto{g'(j+k-2)}{f(i+k-1)} \neq \equalto{g'(j+k-1)}{g(i+k-2)} \neq \equalto{g'(j+k)}{g(i+k)}$$ 

        \item
        [$\begin{pmatrix}
        S & M... & M\\
        M & M... & S\end{pmatrix}$] is analogous to type $\begin{pmatrix}
        M & M... & S\\
        S & M... & M\end{pmatrix}$ if we swap the roles of $f$ and $g$.

        \item
        [$\begin{pmatrix}
        S & M... & S\\
        S & M... & M\end{pmatrix}$] is analogous to type $\begin{pmatrix}
        M & M... & M\\
        S & M... & S\end{pmatrix}$ if we swap the roles of $f$ and $g$.
    \end{itemize}
    \item Type $(P)$ is one of the following two types of movement:
    \begin{itemize}
        \item
        [$\begin{pmatrix}
        M \\
        S \end{pmatrix}$] which means we have
        \begin{tabular}{c}
        $f(i) \neq f(i+1)$ \\
        $g(i) = g(i+1)$ \\
        \end{tabular}, so we define:\\ 

        $$\equaltog{f'(j)}{f(i)} \neq \equaltog{f'(j+1)}{f(i+1)}$$ 
        $$\equalto{g'(j)}{g(i)} \neq \equalto{g'(j+1)}{\text{any neighbor of } g(i+1)}$$ 

        \item
        [$\begin{pmatrix}
        S \\
        M \end{pmatrix}$] is analogous to type $\begin{pmatrix}
        M \\
        S \end{pmatrix}$ if we swap the roles of $f$ and $g$.
    \end{itemize}
    \item Type $(P,A,...,A)$ is one of the following two types of movement:
    \begin{itemize}
        \item
        [$\begin{pmatrix}
        M & M...\\
        S & M...\end{pmatrix}$] which means we have
        \begin{tabular}{c}
        $f(i) \neq f(i+1) \neq ... \neq f(i+k-1) \neq f(i+k)$ \\
        $g(i) = g(i+1) \neq ... \neq g(i+k-1) \neq g(i+k)$ \\
        \end{tabular}, so we define:\\ 

        $$\equaltog{f'(j)}{f(i)} \neq \equaltog{f'(j+1)}{f(i+1)} \neq ... \neq \equaltog{f'(j+k-1)}{f(i+k-1)} \neq \equaltog{f'(j+k)}{f(i+k)}$$ 
        $$\equalto{g'(j)}{g(i)} \neq \equalto{g'(j+1)}{g(i+2)} \neq ... \neq \equalto{g'(j+k-1)}{g(i+k)} \neq \equalto{g'(j+k)}{g(i+k-1)}$$ 

        \item
        [$\begin{pmatrix}
        S & M...\\
        M & M...\end{pmatrix}$] is analogous to type $\begin{pmatrix}
        M & M...\\
        S & M...\end{pmatrix}$ if we swap the roles of $f$ and $g$.
    \end{itemize}
\end{itemize}
Note that for each block, $f'$ and $g'$ visit all of the edges that $f$ and $g$ visit, so complete $f'$ and $g'$ are $k$-sweeps in $G$, for some $k\in\mathbb{N}$. Let us comment on the distance that $f'$ and $g'$ maintain. For blocks of type $(A)$, $(P,P)$ and $(P,A,...,A,P)$ lazy sweeps $f$ and $g$ start and end in the same vertices as sweeps $f'$ and $g'$, ensuring that the next block can start at a distance of at least $d$. The final block may be of type $(P)$ or $(P,A,...,A)$ for which sweeps $f'$ and $g'$ start in the same vertices as lazy sweeps $f$ and $g$, at a distance of at least $d$. They do not end in the same vertices as $f$ and $g$ but since these are the last steps of $f'$ and $g'$, it is enough that their distance is at least $d-1$, which is clear from the construction. So for sweeps $f'$ and $g'$ we have $m_G(f',g')\geq d-1$, which proves $\sigma^{\times}_{E}(G) \geq \sigma^{\boxtimes}_{E}(G)-1$.   
\end{proof}

   Let us observe an example of the construction of two $9$-sweeps $f'$ and $g'$ from two lazy $10$-sweeps $f$ and $g$, by the method described in the previous proof (Figure \ref{primjer}).
\begin{figure}[H]
    \centering
    $$f(1) \neq f(2) \neq f(3) = f(4) \neq f(5) \neq f(6) = f(7) = f(8) = f(9) \neq f(10)$$
    $$g(1) = g(2) \neq g(3) \neq g(4) \neq g(5) \neq g(6) \neq g(7) \neq g(8)\neq g(9) \neq g(10)$$
    \begin{tabular}{|c|c|c|c|c|c|c|c|c|c|}
    \hline
       $f$ & $M$ & $M$ & $S$ & $M$ & $M$ & $S$ & $S$ & $S$ & $M$\\
       $g$ & $S$ & $M$ & $M$ & $M$ & $M$ & $M$ & $M$ & $M$ & $M$\\
    \hline
        & $\downarrow$ & $\downarrow$ & $\downarrow$ & $\downarrow$ & $\downarrow$ & $\downarrow$ & $\downarrow$ & $\downarrow$ & $\downarrow$\\
    \hline
        activity & $P$ & $A$ & $P$ & $A$ & $A$ & $P$ & $P$ & $P$ & $A$\\
   \hline
        \multicolumn{10}{c}{$\downarrow$}\\
    \hline
        blocks & \multicolumn{3}{|c|}{$P A P$} & $A$ & $A$ & \multicolumn{2}{|c|}{$P P$} & \multicolumn{2}{|c|}{$P A$}\\
   \hline    
    \end{tabular}\\

    \vspace{5pt}
    \begin{tabular}{|c|c|}
    \hline
    \multicolumn{2}{|c|}{block $P A P 
    \begin{pmatrix}
        M & M & S\\
        S & M & M
    \end{pmatrix}$}\\
    \hline
    $f'(1)=f(1)$ & $g'(1)=g(1)$\\
    $f'(2)=f(2)$ & $g'(2)=g(3)$\\
    $f'(3)=f(4)$ & $g'(3)=g(4)$\\
    \hline
    \multicolumn{2}{|c|}{block $A \begin{pmatrix}
        M\\
        M
    \end{pmatrix}$}\\
    \hline
    $f'(3)=f(4)$ & $g'(3)=g(4)$\\
    $f'(4)=f(5)$ & $g'(4)=g(5)$\\
    \hline
    \multicolumn{2}{|c|}{block $A \begin{pmatrix}
        M\\
        M
    \end{pmatrix}$}\\
    \hline
    $f'(4)=f(5)$ & $g'(4)=g(5)$\\
    $f'(5)=f(6)$ & $g'(5)=g(6)$\\
    \hline
    \multicolumn{2}{|c|}{block $P P \begin{pmatrix}
        S & S\\
        M & M
    \end{pmatrix}$}\\
    \hline
    $f'(5)=f(6)$ & $g'(5)=g(6)$\\
    $f'(6)=f(5)$ & $g'(6)=g(7)$\\
    $f'(7)=f(8)$ & $g'(7)=g(8)$\\
    \hline
    \multicolumn{2}{|c|}{block $P A \begin{pmatrix}
        S & M\\
        M & M
    \end{pmatrix}$}\\
    \hline
    $f'(7)=f(8)$ & $g'(7)=g(8)$\\
    $f'(8)=f(10)$ & $g'(8)=g(9)$\\
    $f'(9)=f(9)$ & $g'(9)=g(10)$\\
    \hline
    \end{tabular}\\ 
    $$\equaltog{f'(1)}{f(1)} \neq \equaltog{f'(2)}{f(2)} \neq \equaltog{f'(3)}{f(4)} \neq \equaltog{f'(4)}{f(5)} \neq \equaltog{f'(5)}{f(6)} \neq \equaltog{f'(6)}{f(5)} \neq \equaltog{f'(7)}{f(8)} \neq \equaltog{f'(8)}{f(10)} \neq \equaltog{f'(9)}{f(9)}$$
    $$\equalto{g'(1)}{g(1)} \neq \equalto{g'(2)}{g(3)} \neq \equalto{g'(3)}{g(4)} \neq \equalto{g'(4)}{g(5)} \neq \equalto{g'(5)}{g(6)} \neq \equalto{g'(6)}{g(7)} \neq \equalto{g'(7)}{g(8)} \neq \equalto{g'(8)}{g(9)} \neq \equalto{g'(9)}{g(10)}$$ 
    \caption{Example of two lazy $10$-sweeps $f$ and $g$, their movement sequences, derived activity sequence, partition into blocks and  construction of $9$-sweeps $f'$ and $g'$}
    \label{primjer}
\end{figure}

This result also holds for vertex spans, due to every $k$-sweep being a $k$-track, for any $k\in\mathbb{N}$, so we have:

\begin{corollary}\label{kor_vertex}
    Let $G$ be a graph. Then 
     $$\sigma^{\times}_{V}(G) \geq \sigma^{\boxtimes}_{V}(G)-1.$$
\end{corollary}

Finally, we prove a strong relation between all three edge spans.

\begin{theorem}
 Let $G$ be a graph. Then
    $$\max\{\sigma^{\times}_{E}(G),\sigma^{\square}_{E}(G)\}\leq\sigma^{\boxtimes}_{E}(G)\leq \max\{\sigma^{\times}_{E}(G),\sigma^{\square}_{E}(G)\}+1.$$
\end{theorem}
\begin{proof}
    The left inequality is the claim of Proposition \ref{strong_max}. For the right inequality, we join the result from Theorem \ref{analog}, $|\sigma^{\times}_E(G)-\sigma^{\square}_E(G)| \leq 1$ and the result from Theorem \ref{strong_direct}, $\sigma^{\times}_{E}(G) \geq \sigma^{\boxtimes}_{E}(G)-1$. We distinguish three cases:
    \begin{enumerate}
        \item $\sigma^{\times}_E(G)=\sigma^{\square}_E(G)$. In this case $\sigma^{\boxtimes}_{E}(G)\leq \sigma^{\times}_{E}(G)+1=\sigma^{\square}_E(G)+1$, so $\sigma^{\boxtimes}_{E}(G)\leq \max\{\sigma^{\times}_{E}(G),\sigma^{\square}_{E}(G)\}+1$ follows.
        \item $\sigma^{\times}_E(G)=\sigma^{\square}_E(G)+1$. Now $\sigma^{\times}_E(G)=\max\{\sigma^{\times}_{E}(G),\sigma^{\square}_{E}(G)\}$, so the claim follows from Theorem \ref{strong_direct}.
        \item $\sigma^{\times}_E(G)=\sigma^{\square}_E(G)-1$. We have
        $$\sigma^{\boxtimes}_{E}(G)\leq \sigma^{\times}_{E}(G)+1=\sigma^{\square}_E(G)-1+1=\sigma^{\square}_E(G),$$ so again $\sigma^{\boxtimes}_{E}(G)\leq \max\{\sigma^{\times}_{E}(G),\sigma^{\square}_{E}(G)\}+1$ follows.
    \end{enumerate}
\end{proof}

An analogous result holds for vertex spans:

\begin{corollary}
    Let $G$ be a graph. Then
    $$\max\{\sigma^{\times}_{V}(G),\sigma^{\square}_{V}(G)\}\leq\sigma^{\boxtimes}_{V}(G)\leq \max\{\sigma^{\times}_{V}(G),\sigma^{\square}_{V}(G)\}+1.$$
\end{corollary}
\begin{proof}
    The inequality on the left is Proposition 4.1. in \cite{nas}. For the right inequality, we distinguish two cases:
    \begin{enumerate}
        \item $\sigma^{\times}_{V}(G)\geq \sigma^{\square}_{V}(G)$. Now $\max\{\sigma^{\times}_{V}(G),\sigma^{\square}_{V}(G)\}+1=\sigma^{\times}_{V}(G)+1$, so the claim follows from Corollary \ref{kor_vertex}.
        \item $\sigma^{\square}_{V}(G)>\sigma^{\times}_{V}(G)$. Now, from Theorem 4.2. in \cite{nas} we have $\sigma^{\square}_{V}(G)=\sigma^{\times}_{V}(G)+1$, so $\max\{\sigma^{\times}_{V}(G),\sigma^{\square}_{V}(G)\}+1=\sigma^{\times}_{V}(G)+2$. It follows $$\sigma^{\boxtimes}_{V}(G)\leq \sigma^{\times}_{V}(G)+1<\sigma^{\times}_{V}(G)+2=\max\{\sigma^{\times}_{V}(G),\sigma^{\square}_{V}(G)\}+1.$$
    \end{enumerate}
\end{proof}


\section{Conclusion and further work}\label{con}
In this paper, we define edge span variants of a graph analogously to vertex span variants defined in \cite{nas}, both based on the work in \cite{banic}. We prove several relations between different edge spans of a graph, some analogous to the results for vertex spans, and some new and important connections between all three spans that bound the strong edge span from above. Moreover, we also prove the result for vertex spans, improving previously known results. We show that the upper bound is tight, by finding a graph family for which the strong edge span does not equal the larger value between the direct edge span and the Cartesian edge span. We conjecture that this does not hold for vertex spans; we believe that for any graph $G$, $\sigma^{\boxtimes}_{V}(G)= \max\{\sigma^{\times}_{V}(G),\sigma^{\square}_{V}(G)\}$. This problem remains open for future research.

\textbf{Author Contributions:} Conceptualization, A.Š. and T.V.; Investigation, A.Š. and T.V.; Data curation, A.Š. and T.V.; Writing—original draft, A.Š. and T.V.; Writing—review \& editing, A.Š. and T.V. All authors have read and agreed to the published version of the manuscript.\\

\textbf{Funding:} This research received no external funding.\\

\textbf{Data Availability Statement:} Not applicable.\\

\textbf{Acknowledgments:} Not applicable.\\

\textbf{Conflicts of Interest:} The authors declare that they have no conflict of interest.

\begin{thebibliography}{33}

\bibitem{about}Aboutahoun, A., Mahdi, S., El-Alem, M., \& ALrashidi, M., Modified and Improved Algorithm for Finding a Median Path with a Specific Length ($\ell$) for a Tree Network. Mathematics, 11(16), 3585. (2023).


\bibitem{ambrosino2}Ambrosino, D., \& Cerrone, C., A Rich Vehicle Routing Problem for a City Logistics Problem. Mathematics, 10(2), 191. (2022).

\bibitem{aziz} Aziz, F., Wilson, R. C., \& Hancock, E. R., Backtrackless walks on a graph. IEEE transactions on neural networks and learning systems, 24(6), 977-989. (2013).

\bibitem{banic} Banič, I.  \& Taranenko, A., Span of a graph: keeping the safety distance. Discrete Mathematics and Theoretical Computer Science, 25 (Graph Theory). (2023).

\bibitem{dash}Dashtbali, M., Malek, A., \& Mirzaie, M., Optimal control and differential game solutions for social distancing in response to epidemics of infectious diseases on networks. Optimal Control Applications and Methods, 41(6), 2149-2165. (2020).

\bibitem{taranenko1}Dravec, T., Mikalački, M., \& Taranenko, A. (2024). Graphs with span 1 and shortest optimal walks. arXiv preprint arXiv:2410.19524.

\bibitem{nas} Erceg, G., Šubasic, A. \& Vojković, T., Some results on the maximal safety distance in a graph. FILOMAT, 37(15), 5123--5136. (2023).


\bibitem{kinez} Edmonds, J.  \& Johnson, E. L., Matching, Euler tours and the Chinese postman. Mathematical programming, 5, 88--124. (1973).

\bibitem{taranenko2}Grašič, M., Mouron, C., \& Taranenko, A. (2024). The strong vertex span of trees. arXiv preprint arXiv:2412.03266.


\bibitem{lagos}Lagos, A. R., Kordonis, I., \& Papavassilopoulos, G. P., Games of social distancing during an epidemic: local vs statistical information. Computer Methods and Programs in Biomedicine Update, 2, 100068. (2022).


\bibitem{lelek}Lelek, A., Disjoint mappings and the span of spaces. Fund. Math., 55, 199--214. (1964).

\bibitem{patasius}Patašius, M., Šimkienė, J., Sokas, D., \& Pranskūnas, A., Method for Finding the Limits of Blood Vessel Landmarks in Eye Fundus Images Based on Distances in Graphs: Preliminary Results. In XV Mediterranean Conference on Medical and Biological Engineering and Computing–MEDICON 2019: Proceedings of MEDICON 2019, September 26-28, 2019, Coimbra, Portugal (pp. 358-366). Springer International Publishing. (2020).


\bibitem{rumi}Rumiantsev, B. V., Kochkarov, R. A., \& Kochkarov, A. A., Graph-Clustering Method for Construction of the Optimal Movement Trajectory under the Terrain Patrolling. Mathematics, 11(1), 223. (2023).

\bibitem{rupa}Rupapara, V., Narra, M., Gunda, N. K., Gandhi, S., \& Thipparthy, K. R., Maintaining social distancing in pandemic using smartphones with acoustic waves. IEEE Transactions on Computational Social Systems, 9(2), 605-611. (2021).

\bibitem{multil}Šubašić, A., \& Vojković, T. Vertex Spans of Multilayered Cycle and Path Graphs. Axioms, 13.4, (2024): 236.

\bibitem{gross} West D.B., Introduction to graph theory (Vol. 2), Prentice hall, Upper Saddle River, (2001).


\end{thebibliography}
\end{document}